\documentclass[11pt,oneside,onecolumn,reqno]{amsart}

\usepackage{amssymb,latexsym,mathrsfs,
amsxtra,amscd,amsfonts,amsthm,amsmath,verbatim,epsfig}

\usepackage{color}

\usepackage[colorlinks,pagebackref=true]{hyperref}

\usepackage[english]{babel}

\makeatletter

\makeatletter
\def\widebreve#1{\mathop{\vbox{\m@th\ialign{##\crcr\noalign{\kern3\p@}%
      \brevefill\crcr\noalign{\kern3\p@\nointerlineskip}%
      $\hfil\displaystyle{#1}\hfil$\crcr}}}\limits}

\def\brevefill{$\m@th \setbox\z@\hbox{$\braceld$}%
  \bracelu\leaders\vrule \@height\ht\z@ \@depth\z@\hfill\braceru$}
\makeatletter

\def\@citecolor{blue}
\def\@linkcolor{blue}
\def\@urlcolor{blue}
\def\@urlcolor{blue}

\numberwithin{equation}{section}

\def\dim{\operatorname{dim}}

\def\height{\operatorname{ht}}

\def\ann{\operatorname{Ann}}

\def\jr{\operatorname{jr}}
\def\H{\operatorname{H}}
\def\Proj{\operatorname{Proj}}

\def\ZZ{\mathbb Z}
\def\QQ{\mathbb Q}

\newcommand{\NN}{\mathbb N}
\newcommand{\mm}{\mathfrak m}

\newcommand{\ov}{\overline}
\newcommand{\lm}{{\lambda}}

\newcommand{\R}{\mathcal R}

\newcommand{\I}{{\bf I}}
\newcommand{\n}{{\bf{n}}}

\newcommand{\m}{{\bf {m}}}
\newcommand{\rrr}{{\bf{r}}}
\newcommand{\ii}{i=1,\ldots,s}
\newcommand{\idl}{{I_1},\ldots,{I_s}}
\newcommand{\id}{{I_1}\cdots{I_s}}

\newcommand{\lf}{\left(}
\newcommand{\rg}{\right)}
\newcommand{\fil}{\mathcal F}

\newcommand{\ga}{{G(\fil)}}

\newcommand{\po}{P_{\fil}}
\newcommand{\ho}{H_{\fil}}

\newcommand{\bl}{\begin{lemma}}
\newcommand{\el}{\end{lemma}}
\newcommand{\bt}{\begin{theorem}}
\newcommand{\et}{\end{theorem}}
\newcommand{\ben}{\begin{enumerate}}
\newcommand{\een}{\end{enumerate}}
\newcommand{\bpf}{\begin{proof}}
\newcommand{\eepf}{\end{proof}}
\newcommand{\beqn}{\begin{eqnarray*}}
\newcommand{\eeqn}{\end{eqnarray*}}
\newcommand{\beqnn}{\begin{eqnarray}}
\newcommand{\eeqnn}{\end{eqnarray}}
\newcommand{\bd}{\begin{definition}}
\newcommand{\ed}{\end{definition}}
\newcommand{\bp}{\begin{proposition}}
\newcommand{\ep}{\end{proposition}}
\newcommand{\bc}{\begin{corollary}}
\newcommand{\ec}{\end{corollary}}
\newcommand{\bex}{\begin{example}}
\newcommand{\eex}{\end{example}}

\newcommand{\wrt}{with respect to }

\theoremstyle{plain}
\newtheorem{theorem}{Theorem}[section]
\newtheorem{corollary}[theorem]{Corollary}
\newtheorem{proposition}[theorem]{Proposition}
\newtheorem{lemma}[theorem]{Lemma}

\newtheorem{example}[theorem]{Example}
\newtheorem{definition}[theorem]{Definition}

\theoremstyle{remark}
\newtheorem{remark}[theorem]{Remark}
\numberwithin{equation}{theorem}
\begin{document}

\title[Local Cohomology of multi-Rees algebras]{Local cohomology of multi-Rees algebras, Joint reduction numbers  and product of complete ideals}
\author{Parangama Sarkar}
\author{J. K. Verma}
\address{ Mathematics Department, Indian Institute of Technology Bombay, Mumbai, India 400076}

\email{parangama@math.iitb.ac.in}
\email{jkv@math.iitb.ac.in}

\thanks{The first author is supported by CSIR Fellowship of Government of India.} 

\begin{abstract}We find conditions on the local cohomology modules of multi-Rees algebras of admissible filtrations which enable us to predict joint reduction numbers. As a consequence we are able to prove a generalisation of a result of Reid-Roberts-Vitulli
in the setting of analytically unramified local rings for completeness of power products of complete ideals.
\end{abstract}
\maketitle

\thispagestyle{empty}
\section{Introduction}
The objective of this paper is to find suitable conditions on the local cohomology modules of  multi-Rees algebras and associated graded rings of multigraded admissible filtrations  of  ideals in an analytically unramified local ring $(R,\mm)$ and apply these to detect their joint reduction vectors and completeness of products of complete ideals.
\\Recall that if $R$ is a commutative ring and $I$ is an ideal of $R$ then $a \in R$ is called integral over  $I,$ if $a$ is a root of a monic polynomial  $x^n+a_1x^{n-1}+\dots+a_{n-1}x+a_n$  for some $a_i \in I^i$ for $i=1,2,\dots,n.$ The integral closure of $I,$ denoted by $\ov{I},$ is the set of all
elements of $R$ which are integral over $I.$  If $I=\ov{I},$ then $I$ is called complete or integrally closed.
O. Zariski \cite{z} proved that product of complete ideals is complete  in the polynomial ring $k[x,y]$ where $k$ is an algebraically closed  field of characteristic zero. This was generalised
to two-dimensional regular local rings in Appendix $5$ of \cite{zariski}. This result is known as {\em Zariski's Product Theorem.} C. Huneke \cite {h} showed that product of complete ideals is not complete in higher dimensional regular local rings. Since the appearance of this counterexample of Huneke, several results have appeared in the literature which identify  classes of complete ideals in local rings of dimension at least $3$ whose products are complete. The following result  due to L. Reid, L. G. Roberts and M. A. Vitulli \cite[Proposition 3.1]{V} about complete monomial ideals is rather surprising:
\bt
Let $R=k[X_1,\dots,X_d]$ be a polynomial ring of dimension $d \geq 1$ over a field $k.$ Let $I$ be a monomial ideal of $R$  so that  $I^{n}$ is  complete for all $1 \leq n \leq d-1.$ Then $I^{n}$ is  complete for all $n\geq 1.$
\et
This can be thought of as a partial generalisation of the Zariski's  Product Theorem for $d=2.$   This theorem was proved using tools from convex geometry. In this paper, we approach this result using vanishing of  local cohomology modules of multi-Rees algebras and prove   the following result about completeness of power products of $\mm$-primary monomial ideals.
\bt
Let $R=k[X_1,\dots,X_d]$ where $d\geq 1$ and $\mm=(X_1,\dots,X_d)$ be the maximal homogeneous ideal of $R.$ Let $\idl$ be $\mm$-primary monomial ideals of $R.$ Suppose $\I^{\n}$ is complete for all $\n\in\NN^s$ such that $1\leq|\n|\leq d-1.$ Then $\I^{\n}$ is complete for all $\n\in\NN^s$ with $|\n|\geq 1.$ 
\et
We  prove the above result as a consequence of a more general result for complete ideals in analytically unramified  local rings. In order to state this  and other results  proved  in this paper, we recall certain definitions and set up notation.
\\Throughout this paper, $(R,\mm)$ denotes a Noetherian local ring of dimension $d$ with infinite residue field. Let $ I_1,\ldots,I_s$ be $\mm$-primary ideals of $R$ and we denote the collection of these ideals $(\idl)$ by $\I.$ For $s\geq 1,$ we  put ${\bf {e}}=(1,\ldots,1),\; {\bf{0}}=
 (0,\ldots,0)\in{\ZZ}^s$ and for all $\ii,$ ${\bf {e_i}}=(0,\ldots,1,\ldots,0)\in{\ZZ}^s$ where $1$ occurs at $i$th position. 
 For $\n=(n_1,\ldots,n_s)\in{\ZZ}^s,$ we write ${\I}^{\n}=I_{1}^{n_1}\cdots I_{s}^{n_s}$ and $\n^+=(n_1^+,\ldots,n_s^+)$ where $n_i^+=\max\{0,n_i\}.$
  For $s \geq 2$ and  $\alpha=({\alpha}_1,\ldots,{\alpha}_s)\in{\NN}^s,$ we put $|\alpha|={\alpha}_1+\cdots+{\alpha}_s.$ We define $\m=(m_1,\ldots,m_s)\geq\n=(n_1,\ldots,n_s)$ if $m_i\geq n_i$ for all $\ii.$ 
By the phrase ``for all large $\n$", we mean $\n\in\NN^s$ and $n_i\gg 0$ for all $\ii.$
\bd   A set of ideals $\fil=\lbrace\fil(\n)\rbrace_{\n\in \ZZ^s}$ is called a $\ZZ^s$-graded {\bf{$\I$-filtration}} if for all $\m,\n\in\ZZ^s,$
{\rm (i)} ${\I}^{\n}\subseteq\fil(\n),$
 {\rm (ii)} $\fil(\n)\fil(\m)\subseteq\fil(\n+\m)$  and {\rm (iii)} if $\m\geq\n,$ $\fil(\m)\subseteq\fil(\n).$ 
\ed
Let $t_1,\ldots,t_s$ be indeterminates. For $\n\in\ZZ^s,$ we put ${\bf t}^{\n}=t_{1}^{n_1}\cdots t_{s}^{n_s}$ and denote the $\NN^s$-graded {\bf{Rees ring of $\fil$}} by $\mathcal{R}(\fil)=\bigoplus\limits_{\n\in \NN^s}
{\fil}(\n){\bf{t}}^{\n}$ and the $\ZZ^s$-graded {\bf{extended Rees ring of $\fil$}} by 
$\mathcal{R}'(\fil)=\bigoplus \limits_ {\n\in{\ZZ}^s}{\fil}(\n){\bf{t}}^{\n}.$ For an $\NN^s$-graded ring $S=\bigoplus\limits_{\n\geq{\bf{0}}}S_{\n},$ we denote the ideal $\bigoplus \limits_{\n\geq {\bf{e}}}S_{\n}$ by $S_{++}.$ Let $\ga=\bigoplus\limits_{\n\in{\NN}^s}{{\fil}(\n)}/{{\fil}(\n+{\bf e})}$ be the {\bf{associated multigraded ring of $\fil$ with respect to $\fil({\bf e}).$}} For $\fil=\{\I^{\n}\}_{\n \in \ZZ^s}$, we set $\mathcal R(\fil)=\mathcal R(\I),$ $\mathcal R^\prime(\fil)=
\mathcal R^\prime(\I),$ $\ga=G(\I)$ and $\mathcal R(\I)_{++}=\R_{++}.$ 
\bd A $\ZZ^s$-graded $\I$-filtration $\fil=\lbrace\fil(\n)\rbrace_{\n\in \ZZ^s}$ of ideals in $R$ is 
called an $\I$-{\bf{admissible filtration}} if 
${{\fil}(\n)}={\fil}(\n^+)$ for all $\n\in\ZZ^s$ and $\mathcal{R}'(\fil)$
 is a finite $\mathcal{R}'(\I)$
 -module. \ed
The principal examples of admissible filtrations with which we are concerned in this paper   are (i) the  $\I$-adic filtration 
$\{ \I^{\n}\}_{\n \in \ZZ^s}$ in a Noetherian local ring and (ii) the integral closure filtration $\{ \overline{\I^{\n} } \}_{\n \in \ZZ^s}$ in an analytically unramified local ring.  
It is proved in \cite[Proposition 2.5] {msv} that if  $\fil=\lbrace\fil(\n)\rbrace_{\n\in \ZZ^s}$ is an $\I$-admissible filtration of ideals in $R$ then $\R(\fil)$ is a finitely generated $\R(\I)$-module.
\\Recall that an ideal $J$ contained in an ideal $I$ is called a reduction of $I$ if $JI^n=I^{n+1}$ for
all large $n.$ The role that reductions of ideals play in the study of Hilbert-Samuel functions of $\mm$-primary ideals, is played by joint reductions,  introduced by D. Rees in \cite{rees3}, of a sequence of $\mm$-primary ideals $\idl$ to study the multigraded Hilbert-Samuel function $H(\I,\n)=\lm(R/\I^{\n}).$  
 Let ${\bf{q}}=(q_1,q_2,\dots,q_s)\in \NN^s$  and $|{\bf q}|=d \geq 1.$ A set of elements 
$ \{a_{ij}  \in I_i \mid i=1,2,\dots,s; j=1,2,\dots,q_i\}$ is called a joint reduction of  the set of ideals $(I_1,\ldots,I_s)$ of  type ${\bf{q}}$ if
there exists an $\m \in \NN^s$ so that for all  $ \n \geq \m \in \NN^s,$
$$\sum_{i=1}^s\sum_{j=1}^{q_i} a_{ij} 
I_1^{n_1}I_2^{n_2} \dots I_{i-1}^{n_{i-1}}I_i^{n_i-1}I_{i+1}^{n_{i+1}}\dots I_s^{n_s}=
I_1^{n_1}I_2^{n_2}  \dots I_s^{n_s}.$$
The vector  $\m$ is called a joint reduction vector. We estimate joint reduction vectors using local cohomology modules of multi-Rees rings.
In order to achieve this we need to work with joint reductions in a more general setting.
 D. Kirby and Rees \cite{Kirby-Rees} generalised it further in the setting of multigraded rings and modules which we recall next.
\bd Let $R=\bigoplus\limits_{\n \in \NN^s}R_{\n}$ be a standard Noetherian $\NN^s$-graded ring defined over a local ring $(R_{\bf{0}},\mm)$  and
$M=\bigoplus\limits_{\n \in \ZZ^s} M_{\n}$ be a finite $\ZZ^s$-graded $R$-module.
 A {\bf{joint reduction of type ${\bf{q}}$ of $R$ \wrt $M$}} is a set of elements $$\mathcal A_{\bf{q}}(M)=\lbrace a_{ij}\in R_{{\bf {e_i}}}:j=1,\ldots,q_i;\ii\rbrace$$ generating an ideal $J$ of $R$ irrelevant \wrt $M,$ i.e. $(JM)_{\n}=M_{\n}$ for all large $\n.$
\ed
Kirby and Rees \cite{Kirby-Rees} proved existence of joint reduction of type ${\bf{q}}$ of $R$ \wrt $M$ if $\displaystyle |{\bf{q}}|\geq \dim\lf\frac{M}{\mm M}\rg+1$ and the residue field $R_{\bf{0}}/\mm$ is infinite. Here $\dim R$ is defined to be $\max(\height P)$ where $P$ ranges over the relevant prime ideals of $R$ if $R$ is not trivial and $-1$ if $R$ is trivial. For $M,$ $\dim M$ is defined to be $\displaystyle\dim\lf\frac{R}{\ann_R M}\rg.$
\\ Let $(R,\mm)$ be a Noetherian local ring of dimension $d\geq 1,$ $\idl$ be $\mm$-primary ideals of $R$ and $\fil=\lbrace\fil(\n)\rbrace_{\n\in \ZZ^s}$ be a $\ZZ^s$-graded $\I$-admissible filtration of ideals in $R.$ Let ${\bf{q}}=(q_1,\ldots,q_s)\in\NN^s$ such that $|{\bf{q}}|=d.$
\bd
A set of elements $\mathcal A_{{\bf{q}}}(\fil)=\lbrace a_{ij}\in I_{i}:j=1,\ldots,q_i;\ii\rbrace$ is called a {\bf{joint reduction of $\fil$}} of type ${\bf{q}}$ if  the set $\lbrace a_{ij}t_i\in \R(\I)_{{\bf {e_i}}}:j=1,\ldots,q_i;\ii\rbrace$ is a joint reduction of type ${\bf{q}}$ of $\R(\I)$ \wrt $\R(\fil),$ i.e. the following equality holds for all  $\n \geq \m $ for   some $\m\in\NN^s:$
$$\sum\limits_{i=1}^{s}{\sum\limits_{j=1}^{q_i}{a_{ij}\fil(\n-{\bf {e_i}})}}=\fil(\n).$$ 
The vector $\m$ is called a {\bf joint reduction vector} of $\mathcal F$ with respect to the
joint reduction $\mathcal A_{{\bf{q}}}(\fil).$
\ed
Let $I,J$ be $\mm$-primary ideals in a Noetherian local ring $(R,\mm)$ of dimension $d\geq 2,$ $\mu,\lambda\geq 1$ and $\mu+\lambda=d.$ For the bigraded filtration $\fil=\lbrace I^rJ^s\rbrace_{r,s\in \ZZ}$ and a joint reduction $\mathcal A(\mu,\lambda)=\lbrace a_i,b_j,\mid a_i\in I, b_j\in J, 1\leq i\leq\mu,1\leq j\leq\lambda\rbrace$, E. Hyry defined the joint reduction number of $\fil$ \wrt $\mathcal A(\mu,\lambda)$ to be the smallest integer $n$ satisfying $$I^{n+1}J^{n+1}=(a_1,\ldots,a_{\mu})I^nJ^{n+1}+(b_1,\ldots,b_{\lambda})I^{n+1}J^n.$$ We  adapt  this definition to define joint reduction number for multigraded filtrations. 
\bd
Let ${\bf{q}}=(q_1,\ldots,q_s)\in\NN^s$ such that $|{\bf{q}}|=d \geq 1.$ The {\bf{joint reduction number of $\fil$ with respect to a joint reduction}} 
$\mathcal A_{{\bf{q}}}(\fil)=\lbrace a_{ij}\in I_i:j=1,\ldots,q_i;\ii\rbrace$ is the smallest integer $n\in\NN,$ denoted by $\jr_{\mathcal A_{{\bf{q}}}}(\fil),$ such that for all $\n\in\NN^s,$
$$\sum\limits_{i=1}^{s}{\sum\limits_{j=1}^{q_i}{a_{ij}\fil \left(\sum\limits_{k\in A}(n+1){\bf {e_k}}+\n-{\bf {e_i}}\right)}}=\fil \left(\sum\limits_{k\in A}(n+1){\bf {e_k}}+\n \right)\mbox{ where } A=\{i | q_i \neq 0\}.$$ We define the {\bf{joint reduction number of $\fil$ of type ${\bf{q}}$}} to be  
$$\jr_{{\bf{q}}}(\fil)=\min\{\jr_{\mathcal A_{{\bf{q}}}}(\fil)\mid\mbox{}\mathcal A_{{\bf{q}}}(\fil)\mbox{ is a joint reduction of }\fil\mbox{ of type }{\bf{q}}\}.$$
\ed
A crucial step in our investigations is to establish a connection between joint reduction vectors and  vanishing of multigraded components of local cohomology modules of multi-Rees algebras. The following result of Hyry plays a crucial role. 
\bl\cite[Lemma 2.3]{hyry}\label{www} Let $S$ be a Noetherian $\ZZ$-graded ring defined over a local ring $(R,\mm).$ 
Let $\mathcal M$ be the homogeneous maximal ideal of $S.$ Let $\mathfrak a\subset \mm$ be an ideal. 
Let $M$ be a finitely generated $\ZZ$-graded $S$-module and $n_0\in\ZZ.$ Then 
$[H^i_\mathcal M(M)]_n=0$ for all $n\geq n_0$ and $i\geq 0$ if and only if 
$[H^i_{(\mathfrak a,S_+)} (M)]_n=0$ for all $n\geq n_0$ and $i\geq 0.$
\el
For convenience, inspired by the above result, we introduce an invariant of local cohomology modules of multigraded modules over multigraded rings. Let $R=\bigoplus\limits_{\n\in\NN^s}R_{\n}$ be a standard Noetherian $\NN^s$-graded ring defined over a local ring $(R_{{\bf{0}}},\mm).$  Let $M=\bigoplus\limits_{\n\in\ZZ^s}M_{\n}$ be a finitely generated $\ZZ^s$-graded $R$-module.
 \bd
 Let $\m \in \ZZ^s.$ We say that the module $M$ satisfies {\bf{Hyry's condition $\H_R(M,\m)$}} if 
 $$[H_{R_{++}}^i(M)]_{\n}=0 \text{ for all } i\geq 0 \text{ and }  \n \geq \m .$$
 \ed
 Suppose $R_{{\bf {e_i}}}\neq 0$ for all $\ii.$ Let $\mathcal M$ be the maximal homogeneous ideal of $R,$ for each $\ii,$ $\mathcal M_i$ be the ideal of $R$ generated by $R_{{\bf {e_i}}}$ and $R_{++}=\bigcap\limits_{i}\mathcal M_i.$
\\Let $I$ be any subset of $\{1,\dots,s\}$ and $J$ be a non-empty subset of $\{1,\dots,s\}.$ Then for  disjoint $I$ and $J,$ we define $\mathcal M_{I,J}=(\bigcap\limits_{i\in I}\mathcal M_i)\bigcap(\sum\limits_{j\in J}\mathcal M_j).$ We prove a multigraded version of the above result due to Hyry.
\bp
Let $R=\bigoplus\limits_{\n\in\NN^s}R_{\n}$ be a standard $\NN^s$-graded Noetherian ring defined over a local ring $(R_{{\bf{0}}},\mm),$ $R_{{\bf {e_i}}}\neq 0$ for all $\ii$ and $M=\bigoplus\limits_{\n\in\NN^s}M_{\n}$ be a finitely generated $\NN^s$-graded $R$-module. Let $I$ be any subset of $\{1,\dots,s\}$ and $J$ be a non-empty subset of $\{1,\dots,s\}$ such that  $I$ and $J$ are disjoint. Suppose ${\bf{a}}=(a_1,\dots,a_s)\in\ZZ^s$ and $[H^{i}_{\mathcal M}(M)]_{\n}=0$ for all $i\geq 0$ and $\n\in\ZZ^s$ such that $n_k>a_k$ for at least one $k\in\{1,\dots,s\}.$ Then $[H^{i}_{\mathcal M_{I,J}}(M)]_{\n}=0$ for all $i\geq 0$ and $\n\geq\bf{ a+e}.$ In particular, $M$ satisfies Hyry's condition $\H_R(M,\bf{ a+e}).$
\ep
 
In order to detect joint reduction vectors of multigraded admissible filtrations, we  use the theory of filter-regular sequences for multigraded modules.

\bd A homogeneous element $a\in R$ is called an {\bf{$M$-filter-regular}} if $(0:_{M}a)_{\n}=0$ for all large $\n.$
Let $a_1,\dots,a_r\in R$ be homogeneous elements. Then $a_1,\dots,a_r$ is called an {\bf{$M$-filter-regular sequence}} if $a_i$ is $M/(a_1,\dots,a_{i-1})M$-filter-regular for all $i=1,\dots,r.$
\ed

\bt
Let $(R,\mm)$ be a Noetherian local ring of dimension $d\geq 1$ and $\idl$ be $\mm$-primary ideals in $R.$ Let $\mathcal{F}=\lbrace\fil(\n)\rbrace_{\n\in{\ZZ}^s}$ be an $\I$-admissible filtration of ideals in $R.$ Suppose $G(\fil)$ satisfies
Hyry's condition $\H_{G(\I)}(G(\fil), \m).$ 
Let ${\bf{q}}\in\NN^s$ such that $|{\bf{q}}|=d$ and $\lbrace a_{ij}\in I_i:j=1,\ldots,q_i;\ii\rbrace$ be a joint reduction of $\fil$ of type ${\bf{q}}$ such that $a^*_{11},\dots,a^*_{1q_1},\dots,a^*_{s1},\dots,a^*_{sq_s}$ is a $G(\fil)$-filter-regular sequence where $a^*_{ij}$ is the image of $a_{ij}$ in $G(\I)_{{\bf {e_i}}}$ for all $j=1,\dots,q_i$ and $\ii.$ Then $$\fil(\n)=\sum_{i=1}^{s}{\sum_{j=1}^{q_i}{a_{ij}\fil(\n-{\bf{{\bf {e_i}}}})}}\mbox{ for all }\n\geq \m+{\bf{q}}.$$
\et
We can now state the main theorem of this paper which gives a generalisation of Reid-Roberts-Vitulli Theorem for zero-dimensional monomial ideals.
\bt
Let $(R,\mm)$ be an analytically unramified Noetherian local ring of dimension $d\geq 2$  and $\idl$ be $\mm$-primary ideals in $R.$ Let $\overline{\R}(\I)=\bigoplus\limits_{\n\in\NN^s}\overline{\I^{\n}}$ satisfy the condition  
$\H_{\R{(\I)}}(\overline{\R}(\I),{\bf{0}}).$ 
Suppose $\I^{\n}$ is complete for all $\n\in\NN^s$ such that $1\leq|\n|\leq d-1.$ Then $\I^{\n}$ is  complete for all $\n\in\NN^s$ with  $|\n|\geq 1.$ 
\et

We prove that if $\overline{\R}(\I)$ is Cohen-Macaulay then it satisfies Hyry's condition   $\H_{\R{(\I)}}(\overline{\R}(\I),{\bf{0}}).$ By Hochster's Theorem \cite[Theorem 6.3.5]{BH} about Cohen-Macaulayness of normal semigroup rings,  $\overline{\R}(\I)$ is Cohen-Macaulay if $\I$ consists of monomial ideals in a polynomial ring over a field.

\section{Existence of joint reductions consisting of filter-regular sequences}
Let $R=\bigoplus\limits_{\n\in\NN^s}R_{\n}$ be a standard Noetherian $\NN^s$-graded ring defined over an Artinian local ring $(R_{{\bf{0}}},\mm).$ Let $M=\bigoplus\limits_{\n\in\NN^s}M_{\n}$ be a finitely generated $\NN^s$-graded $R$-module. Let ${\Proj}^s(R)$ denote the set of all homogeneous prime ideals $P$ in $R$ such that $R_{++}\nsubseteq P$ and $M^{\Delta}=\bigoplus\limits_{n\geq 0}M_{n{\bf{e}}}.$ By \cite[Theorem 4.1]{hhrz}, there exists a numerical polynomial $P_M\in \QQ[X_1,\dots,X_s]$ of total degree $\dim M^{\Delta}-1$ of the form $$P_M(\n)=\displaystyle\sum\limits_{\substack{\alpha=({\alpha}_1,\ldots,{\alpha}_s)\in{\NN}^s \\ 
|\alpha|\leq \dim M^{\Delta}-1}}(-1)^{\dim M^{\Delta}-1-|{\alpha}|}{e_{\alpha}}(M)\binom{{n_1}+{{\alpha}_1}-1}{{\alpha}_1}
\cdots\binom{{n_s}+{{\alpha}_s}-1}{{\alpha}_s}$$ such that $e_{\alpha}(M)\in\ZZ,$ $P_M(\n)=\lm_{R_{{\bf{0}}}}\lf M_{\n}\rg$ for all large $\n$ and $e_{\alpha}(M)\geq 0$ for all $\alpha\in\NN^s$ such that $|\alpha|=\dim M^{\Delta}-1.$
 \bp\label{existence}
Let $R=\bigoplus\limits_{\n\in\NN^s}R_{\n}$ be a standard Noetherian $\NN^s$-graded ring defined over a local ring $(R_{{\bf{0}}},\mm)$ with infinite residue field $R_{{\bf{0}}}/\mm.$ Let $M=\bigoplus\limits_{\n\in\NN^s}M_{\n}$ be a finitely generated $\NN^s$-graded $R$-module and $\dim M^{\Delta}\geq 1.$ Fix $i\in\{1,\ldots,s\}.$ If $R_{\bf e_i}\neq 0$ then there exists $a\in R_{\bf e_i}$ such that $a$ is $M$-filter-regular. 
\ep
\bpf
Denote $\displaystyle{{M}/{H^{0}_{R_{++}}(M)}}$ by $M'.$ Then $Ass(M')=Ass(M)\setminus V(R_{++}).$ Let $Ass(M')=\{P_1,\ldots,P_k\}.$
\\Let ${\mathfrak a}_i$ be the ideal of $R$ generated by $R_{\bf e_i}.$  Therefore for all $j=1,\ldots,k,$ $P_j\nsupseteq{\mathfrak a}_i.$ Consider the $R_{{\bf{0}}}/\mm$-vector space $R_{\bf e_i}/\mm R_{\bf e_i}.$ Then for each $j=1,\ldots,k,$ $$(P_j\cap R_{\bf e_i}+\mm R_{\bf e_i})/\mm R_{\bf e_i}\neq R_{\bf e_i}/\mm R_{\bf e_i}.$$ Since $R_{{\bf{0}}}/\mm$ is infinite, there exists $a\in R_{{\bf{e_i}}}\setminus \bigcup_{j=1}^k(P_j\cap R_{\bf e_i}+\mm R_{\bf e_i}).$
\\By \cite[Proposition 4.1]{msv}, there exists $\m$ such that $[{H^{0}_{R_{++}}(M)}]_{\n}=0$ for all $\n\geq\m.$ Let $\n\geq\m$ and $x\in (0:_{M}a)_{\n}.$ Then $ax'=0$ in $M'$ where $x'$ is the image of $x$ in $M'.$ Since $a$ is a nonzerodivisor of $M',$ $x\in [{H^{0}_{R_{++}}(M)}]_{\n}=0.$  
\eepf
\bp\label{length}
Let $R=\bigoplus\limits_{\n\in\NN^s}R_{\n}$ be a standard Noetherian $\NN^s$-graded ring defined over an Artinian local ring $(R_{{\bf{0}}},\mm).$ Let $M=\bigoplus\limits_{\n\in\NN^s}M_{\n}$ be a finitely generated $\NN^s$-graded $R$-module. Let $a_i\in R_{{\bf {e_i}}}$ be an $M$-filter-regular element. Then for all large $\n,$ $$\lm_{R_{{\bf{0}}}}\lf\frac{M_{\n}}{a_iM_{\n-{\bf {e_i}}}}\rg=\lm_{R_{{\bf{0}}}}\lf M_{\n}\rg-\lm_{R_{{\bf{0}}}}\lf M_{\n-{\bf {e_i}}}\rg$$ and hence for all $\n\in\ZZ^s,$ $P_{{M}/{a_iM}}(\n)=P_M(\n)-P_M(\n-{\bf {e_i}}).$
\ep
\bpf
Consider the exact sequence of $R$-modules $$0\longrightarrow (0:_{M}a_i )_{\n-{\bf {e_i}}}\longrightarrow M_{\n-{\bf {e_i}}}\overset{a_i}\longrightarrow M_{\n}\longrightarrow \frac{M_{\n}}{a_iM_{\n-{\bf {e_i}}}}\longrightarrow 0.$$ Since $a_i$ is $M$-filter-regular, for all large $\n,$ we get $$\lm_{R_{{\bf{0}}}}\lf\frac{M_{\n}}{a_iM_{\n-{\bf {e_i}}}}\rg=\lm_{R_{{\bf{0}}}}\lf M_{\n}\rg-\lm_{R_{{\bf{0}}}}\lf M_{\n-{\bf {e_i}}}\rg$$ and hence for all $\n\in\ZZ^s,$ $P_{{M}/{a_iM}}(\n)=P_M(\n)-P_M(\n-{\bf {e_i}}).$
\eepf
\bt\label{filter}
Let $R=\bigoplus\limits_{\n\in\NN^s}R_{\n}$ be a standard Noetherian $\NN^s$-graded ring defined over an Artinian local ring $(R_{{\bf{0}}},\mm)$ with infinite residue field $R_{{\bf{0}}}/\mm$ and $R_{\bf e_i}\neq 0$ for all $\ii.$ Let $M=\bigoplus\limits_{\n\in\NN^s}M_{\n}$ be a finitely generated $\NN^s$-graded $R$-module and $\dim M^{\Delta}\geq 1.$ Let $e_{\alpha}(M)>0$ for all $\alpha\in\NN^s$ such that $|\alpha|=\dim M^{\Delta}-1.$ Then for any ${\bf{q}}=(q_1,\dots,q_s)\in\NN^s$ such that $|{\bf{q}}|=\dim M^{\Delta},$ there exist $a_{i1},\dots,a_{iq_i}\in R_{{\bf{e_i}}}$  for all $\ii,$ such that 
 $a_{11},\dots,a_{1q_1},\dots,a_{s1},\dots,a_{sq_s}$ is an $M$-filter-regular sequence and for all large $\n,$ $ M_{\n}=\sum\limits_{i=1}^{s}{\sum\limits_{j=1}^{q_i}{a_{ij}M_{\n-{\bf{e_i}}}}}.$
 \et
\bpf
We use induction on $\dim M^{\Delta}=l.$ Let $l=1.$ Then by Proposition \ref{existence}, for each $\ii,$ there exists $a_i\in R_{{\bf {e_i}}}$ such that $a_i$ is $M$-filter-regular. Since $l=1,$ $P_M(\n)$ is polynomial of total degree zero. Therefore by Proposition \ref{length}, $\displaystyle\lm_{R_{{\bf{0}}}}\lf{M_{\n}}/{a_iM_{\n-{\bf{e_i}}}}\rg=0$ for all large $\n$ and hence we get the required result. Suppose $l\geq 2$ and the result is true for all finitely generated $\NN^s$-graded $R$-module $T$ such that $1\leq\dim T^{\Delta}\leq l-1$ and $e_{\alpha}(T)>0$ for all $\alpha\in\NN^s$ such that $|\alpha|=\dim T^{\Delta}-1.$ Let $M$ be  finitely generated $\NN^s$-graded $R$-module such that $\dim M^{\Delta}=l$ and $e_{\alpha}(M)>0$ for all $\alpha\in\NN^s$ such that $|\alpha|=l-1.$ Fix ${\bf{q}}=(q_1,\dots,q_s)\in\NN^s$ such that $|{\bf{q}}|=\dim M^{\Delta}.$ Let $i=\min\{j\mid q_j\neq 0\}.$ By Proposition \ref{existence}, there exists $a_{i1}\in R_{\bf e_i}$ such that $a_{i1}$ is an $M$-filter-regular element. Let $\displaystyle N={M}/{a_{i1}M}.$ Since $e_{\alpha}(M)> 0$ for all ${\alpha}\in\NN^s$ such that $|\alpha|=\dim M^{\Delta}-1,$ by Proposition \ref{length}, $P_N(\n)$ is a polynomial of degree $l-2$ and hence $\dim N^{\Delta}=l-1.$ Let $\beta=(\beta_1,\dots,\beta_s)\in\NN^s$ such that $|\beta|=l-2$ and $\alpha=\beta+{\bf {e_i}}.$ Then $e_{\beta}(N)=e_{\alpha}(M)>0.$ Let $\m={\bf{q-e_i}}\in\NN^s,$ i.e. $m_i=q_i-1$ and for all $j\neq i,$ $m_j=q_j.$ Since $|\m|=l-1,$ by induction hypothesis there exist $b_{j1},\dots,b_{jm_j}\in R_{{\bf {e_j}}}$  for all $j=1,\dots,s$ such that $b_{11},\dots,b_{1m_1},\dots,b_{s1},\dots,b_{sm_s}$ is an $N$-filter-regular sequence and for all large $\n,$ $\displaystyle N_{\n}=\sum\limits_{k=1}^{s}{\sum\limits_{j=1}^{m_k}{b_{kj}N_{\n-{\bf{e_k}}}}}.$ Let $a_{ik}=b_{i(k-1)}$ for all $k=2,\ldots,,q_i$ and for all $j\neq i,$ $a_{jk}=b_{jk}$ for all $k=1,\ldots,,q_j.$ Then for all large $\n,$ $\displaystyle M_{\n}=\sum\limits_{i=1}^{s}{\sum\limits_{j=1}^{q_i}{a_{ij}M_{\n-{\bf{e_i}}}}}.$
Since $a_{i1}$ is $M$-filter-regular, $a_{11},\dots,a_{1q_1},\dots,a_{s1},\dots,a_{sq_s}$ is an $M$-filter-regular sequence. \eepf

\bt\label{joint}
Let $(R,\mm)$ be a Noetherian local ring of dimension $d\geq 1$ and $\idl$ be $\mm$-primary ideals in $R.$ Let $\mathcal{F}=\lbrace\fil(\n)\rbrace_{\n\in{\ZZ}^s}$ be an $\I$-admissible filtration of ideals in $R$ and ${\bf{q}}=(q_1,\dots,q_s)\in\NN^s$ such that $|{\bf{q}}|=d.$ Then there exists a joint reduction $\lbrace a_{ij}\in I_i:j=1,\ldots,q_i;\ii\rbrace$ of $\fil$ of type ${\bf q}$ such that
$a^*_{11},\dots,a^*_{1q_1},\dots,a^*_{s1},\dots,a^*_{sq_s}$ is a $G(\fil)$-filter-regular sequence where $a^*_{ij}$ is the image of $a_{ij}$ in $G(\I)_{\bf {e_i}}$ for all $j=1,\ldots,q_i$ and $\ii.$ 
\et

\bpf
Since $\fil$ is an $\I$-admissible filtration, $G(\fil)$ is finitely generated $G(\I)$-module. By \cite[Theorem 2.4]{rees3}, there exists a polynomial $$ \po(\n)=\displaystyle\sum\limits_{\substack{\alpha=({\alpha}_1,\ldots,{\alpha}_s)\in{\NN}^s \\ 
|\alpha|\leq d}}(-1)^{d-|{\alpha}|}{e_{\alpha}}(\fil)\binom{{n_1}+{{\alpha}_1}-1}{{\alpha}_1}
\cdots\binom{{n_s}+{{\alpha}_s}-1}{{\alpha}_s}$$ such that for all large $\n,$ $\po(\n)=\lm_{R}\lf{R}/{\fil(\n)}\rg,$ $e_{\alpha}(\fil)\in\ZZ$ and $e_{\alpha}(\fil)>0$ for all $\alpha\in\NN^s$ where $|\alpha|=d.$ Hence $$\lm\lf\frac{\fil(\n)}{\fil(\n+{\bf e})}\rg=\lm_{R}\lf\frac{R}{\fil(\n+{\bf e})}\rg-\lm_{R}\lf\frac{R}{\fil(\n)}\rg$$ is a numerical polynomial in $\QQ[X_1,\dots,X_s]$ of total degree $d-1$ for all large $\n$ and $e_{\beta}(G(\fil))> 0$ for all $\beta\in\NN^s$ where $|\beta|=d-1.$ Therefore by Theorem \ref{filter}, there exist $a_{i1},\dots,a_{iq_i}\in I_i$  for all $\ii,$ such that 
$a^*_{11},\dots,a^*_{1q_1},\dots,a^*_{s1},\dots,a^*_{sq_s}$ is a $G(\fil)$-filter-regular sequence where $a^*_{ij}=a_{ij}+{\I}^{\bf{e+e_i}}\in G(\I)_{{\bf {e_i}}}$ for all $j=1,\ldots,q_i,$ $\ii$ and $\displaystyle G(\fil)_{\n}=\sum\limits_{i=1}^{s}{\sum\limits_{j=1}^{q_i}{a_{ij}^*G(\fil)_{\n-{\bf {e_i}}}}}$ for all large $\n.$ Hence $$\fil(\n)=\sum\limits_{i=1}^{s}{\sum\limits_{j=1}^{q_i}{a_{ij}\fil(\n-{\bf {e_i}})}}+\fil(\n+{\bf e})\mbox{ for all large }\n.$$ Since $\fil$ is an $\I$-admissible filtration, by \cite{rees3},  for each $\ii,$ there exist an integer $r_i$ such that 
 for all $\n\in\ZZ^s,$ where $n_i\geq r_i,$ $\fil(\n+{\bf {e_i}})=I_i\fil(\n),$
 Hence for all large $\n,$ we get $$\fil(\n)=\sum_{i=1}^{s}{\sum_{j=1}^{q_i}{a_{ij}\fil(\n-{\bf {e_i}})}}+\fil(\n+{\bf e})\subseteq \sum_{i=1}^{s}{\sum_{j=1}^{q_i}{a_{ij}\fil(\n-{\bf {e_i}})}}+\id\fil(\n).$$ Thus by Nakayama's Lemma, for all large $\n,$ $\fil(\n)=\sum\limits_{i=1}^{s}{\sum\limits_{j=1}^{q_i}{a_{ij}\fil(\n-{\bf {e_i}})}}.$
\eepf

\section{Vanishing of Local cohomology modules of Rees algebra of multigraded filtrations}

Let $R=\bigoplus\limits_{\n\in\NN^s}R_{\n}$ be a standard $\NN^s$-graded Noetherian ring defined over a local ring $(R_{{\bf{0}}},\mm)$ and $R_{{\bf {e_i}}}\neq 0$ for all $\ii.$ For a non-empty subset $J$ of $\{1,\dots,s\},$ we define $\mathcal M_J=\sum\limits_{j\in J}\mathcal M_j.$
\bl \label{non-empty}
Let $R=\bigoplus\limits_{\n\in\NN^s}R_{\n}$ be a standard $\NN^s$-graded Noetherian ring defined over a local ring $(R_{{\bf{0}}},\mm),$ $R_{{\bf {e_i}}}\neq 0$ for all $\ii$ and $M=\bigoplus\limits_{\n\in\NN^s}M_{\n}$ be a finitely generated $\NN^s$-graded $R$-module. Let ${\bf{a}}=(a_1,\dots,a_s)\in\ZZ^s$ and $J$ be any non-empty subset of $\{1,\dots,s\}.$ Suppose $[H^{i}_{\mathcal M}(M)]_{\n}=0$ for all $i\geq 0$ and $\n\in\ZZ^s$ such that $n_k>a_k$ for at least one $k\in J.$ Then $[H^{i}_{\mathcal M_J}(M)]_{\n}=0$ for all $i\geq 0$ and $\n\in\ZZ^s$ such that $\sum\limits_{j\in J}n_j>\sum\limits_{j\in J}a_j.$
\el
\bpf
Consider a group homomorphism $\phi:\ZZ^s\rightarrow\ZZ$ defined by $\phi(\n)=\sum\limits_{j\in J}n_j.$ Then $\displaystyle R^{\phi}=\bigoplus\limits_{n\geq 0}(\bigoplus\limits_{\phi(\n)=n}R_{\n}).$ Let $S=(R^{\phi})_{0}$ and $\mathcal N$ be the maximal homogeneous ideal of $S.$ Therefore $(R^{\phi})_{\mathcal N}$ is an $\NN$-graded ring defined over the local ring $S_{\mathcal N}$ and $((\mathcal M_{J})^{\phi})_{\mathcal N}$ is the irrelevant ideal of $(R^{\phi})_{\mathcal N}.$ Then for all $i\geq 0$ and $m>\sum\limits_{j\in J}a_j,$ 
$$
[H^{i}_{(\mathcal M^{\phi})_{\mathcal N}}(M^{\phi})_{\mathcal N}]_{m}= \lf[H^{i}_{\mathcal M^\phi}(M^{\phi})]_{m}\rg_{\mathcal N}=\lf\bigoplus\limits_{\phi(\n)=m}[H^{i}_{\mathcal M}(M)]_{\n}\rg\otimes_{S}{S_{\mathcal N}}.
$$ 
Since $m>\sum\limits_{j\in J}a_j,$ $\phi(\n)=m$ implies $n_k>a_k$ for at least one $k\in J.$ Hence $[H^{i}_{(\mathcal M^{\phi})_{\mathcal N}}(M^{\phi}_{\mathcal N})]_{m}=0$ for all $i\geq 0$ and $m>\sum\limits_{j\in J}a_j.$ Then by Lemma \ref{www}, taking $\mathfrak a=0$ we get $[H^{i}_{((\mathcal M_{J})^{\phi})_{\mathcal N}}(M^{\phi})_{\mathcal N}]_{m}=0$ for all $i\geq 0$ and $m>\sum\limits_{j\in J}a_j.$ Thus $\lf\bigoplus\limits_{\phi(\n)=m}[H^{i}_{\mathcal M_J}(M)]_{\n}\rg\otimes_{S}{S_{\mathcal N}}=0$ for all $i\geq 0$ and $m>\sum\limits_{j\in J}a_j.$ Therefore $[H^{i}_{\mathcal M_J}(M)]_{\n}=0$ for all $i\geq 0$
and $\n\in\ZZ^s$ such that $\sum\limits_{j\in J}n_j>\sum\limits_{j\in J}a_j.$ 
\eepf
\bp\label{multigraded}
Let $R=\bigoplus\limits_{\n\in\NN^s}R_{\n}$ be a standard $\NN^s$-graded Noetherian ring defined over a local ring $(R_{{\bf{0}}},\mm),$ $R_{{\bf {e_i}}}\neq 0$ for all $\ii$ and $M=\bigoplus\limits_{\n\in\NN^s}M_{\n}$ be a finitely generated $\NN^s$-graded $R$-module. Let $I$ be any subset of $\{1,\dots,s\}$ and $J$ be a non-empty subset of $\{1,\dots,s\}$ such that  $I$ and $J$ are disjoint. Suppose ${\bf{a}}=(a_1,\dots,a_s)\in\ZZ^s$ and $[H^{i}_{\mathcal M}(M)]_{\n}=0$ for all $i\geq 0$ and $\n\in\ZZ^s$ such that $n_k>a_k$ for at least one $k\in\{1,\dots,s\}.$ Then $[H^{i}_{\mathcal M_{I,J}}(M)]_{\n}=0$ for all $i\geq 0$ and $\n\geq\bf{ a+e}.$ In particular, $M$ satisfies Hyry's condition $\H_R(M,\bf{ a+e}).$
\ep
\bpf 
We follow the argument given in \cite[Theorem 3.2.6]{ha} and use induction on $r=|I\cup J|.$ Suppose $r=1.$ Since $I,J$ are disjoint and $|J|\geq 1,$ we have $I=\emptyset$ and the result follows from Lemma \ref{non-empty}. Suppose $r\geq 2$ and the result is true upto $r-1.$ Let $I=\{i_1,\dots,i_k\}$ and $J=\{i_{k+1},\dots,i_r\}.$ we use induction on $k.$ If $k=0$ then again by Lemma \ref{non-empty}, we get the result. Suppose $k\geq 1$ and the result is true upto $k-1.$ Let $\mathcal I=I\setminus\{i_k\}$ and $\mathcal J=J\cup\{i_k\}.$ Then $\mathcal M_{\mathcal I,J}+\mathcal M_{\mathcal I,\{i_k\}}=\mathcal M_{\mathcal I,\mathcal J}$ and $\mathcal M_{\mathcal I,J}\cap\mathcal M_{\mathcal I,\{i_k\}}=\mathcal M_{I,J}.$ Consider the following Mayer-Vietoris sequence of local cohomology modules 
$$\displaystyle\cdots\longrightarrow H^{i}_{\mathcal M_{\mathcal I,J}}(M)\bigoplus H^{i}_{\mathcal M_{\mathcal I,\{i_k\}}}(M)\longrightarrow H^{i}_{\mathcal M_{I,J}}(M)\longrightarrow H^{i+1}_{\mathcal M_{\mathcal I,\mathcal J}}(M)\longrightarrow\cdots.$$
Using induction on $k,$ we get $[H^{i+1}_{\mathcal M_{\mathcal I,\mathcal J}}(M)]_{\n}=0$ for all $i\geq 0$ and $\n\geq\bf{ a+e}$ and using induction on $r,$ we get $[H^{i}_{\mathcal M_{\mathcal I,J}}(M)]_{\n}=0=[H^{i}_{\mathcal M_{\mathcal I,\{i_k\}}}(M)]_{\n}$ for all $i\geq 0$ and $\n\geq\bf{ a+e}.$ For $R_{++},$ we take $I=\{1,\dots,s-1\}$ and $J=\{s\}.$
\eepf
Let $R=\bigoplus\limits_{\n\in\NN^s}R_{\n}$ be a standard $\NN^s$-graded Noetherian ring defined over a local ring $(R_{\bf 0},\mm)$ and  let $\mathcal M$ be the maximal homogeneous ideal of $R.$ Let
 $M=\bigoplus\limits_{\n\in\ZZ^s}M_{\n}$ be a finitely generated $\ZZ^s$-graded $R$-module. For all $\ii,$ define \cite{hyry} the $a$-invariants of $M$ as $$ a^i(M)=\sup\{k\in \mathbb{Z}\mid [H^{\dim M}_{\mathcal{M}}(M)]_{\n} \neq 0 \mbox{ for some }\n \in \ZZ^s\mbox{ with }n_i=k \}.$$ 
 Put $a(M)=(a^1(M),\dots,a^s(M)).$
We recall the following result.
\bp\cite[Lemma 3.7]{masuti-tony-verma}\label{minus}
 Let $(R,\mm)$ be a Cohen-Macaulay local ring of dimension $d$, $\idl$ be $\mm$-primary ideals of 
$R$ and $\fil=\lbrace\fil(\n)\rbrace_{\n\in{\ZZ}^s}$ be an $\I$-admissible filtration of ideals in $R.$ Then $a(\R(\fil))=-{\bf{e}}.$
\ep
\begin{remark}\label{property}
Let $(R,\mm)$ be a Noetherian local ring of dimension $d$ and $\idl$ be $\mm$-primary ideals in $R.$ Let $\mathcal{F}=\lbrace\fil(\n)\rbrace_{\n\in{\ZZ}^s}$ be an $\I$-admissible filtration of ideals in $R$ and $\R(\fil)$ be Cohen-Macaulay. Then by Propositions \ref{minus} and \ref{multigraded}, $\R(\fil)$ satisfies Hyry's condition  $\H_{\R{(\I)}}(\R(\fil),{\bf 0}).$ 
\end{remark}
The next theorem is a generalisation of a result due to  Hyry \cite[Theorem 6.1]{hyry1} for $\ZZ^s$-graded admissible filtration of ideals.
\bt
Let $(R,\mm)$ be a Noetherian local ring of dimension $d$ and $\idl$ be $\mm$-primary ideals in $R.$ Let $\mathcal{F}=\lbrace\fil(\n)\rbrace_{\n\in{\ZZ}^s}$ be an $\I$-admissible filtration of ideals in $R$ and let $\R(\fil)$ satisfy Hyry's condition $\H_{\R{(\I)}}(\R(\fil),{\bf 0}).$ Then
\ben
\item $\po(\n)=\ho(\n)$ for all $\n\in\NN^s,$ 
\item for all $\alpha=({\alpha}_1,\dots,{\alpha}_s)\in\NN^s$ such that $|\alpha|=d,$ $$e_{\alpha}(\fil)=\sum\limits_{n_1=0}^{{\alpha}_1}\dots\sum\limits_{n_s=0}^{{\alpha}_s}\binom{{\alpha}_1}{n_1}\dots\binom{{\alpha}_s}{n_s}(-1)^{d-n_1-\dots-n_s}{\lm\lf\frac{R}{\fil(\n)}\rg}$$ where $\n=(n_1,\dots,n_s).$ 
\een
\et
\bpf
$(1)$ Since $\R(\fil)$ satisfies Hyry's condition $\H_{\R{(\I)}}(\R(\fil),{\bf 0}),$  by \cite[Theorem 4.3]{msv}, we get $\po(\n)=\ho(\n)$ for all $\n\in\NN^s.$ 
\\$(2)$ Consider the operators $(\Delta_i\po)(\n)=\po(\n+{\bf e_i})-\po(\n)$ for all $\ii.$ Then $({\Delta}_{s}^{\alpha_s}\cdots{\Delta}_{1}^{\alpha_1}\po)({\bf 0})=e_{\alpha}(\fil)$ for $\alpha=({\alpha}_1,\dots,{\alpha}_s)\in\NN^s,$ $|\alpha|=d.$ By \cite[Proposition 1.2]{Snapper}, $$({\Delta}_{s}^{\alpha_s}\cdots{\Delta}_{1}^{\alpha_1}\po)({\bf 0})=\sum\limits_{n_1=0}^{{\alpha}_1}\dots\sum\limits_{n_s=0}^{{\alpha}_s}\binom{{\alpha}_1}{n_1}\dots\binom{{\alpha}_s}{n_s}(-1)^{d-n_1-\dots-n_s}\po(\n)$$ where $\n=(n_1,\dots,n_s).$ Thus from part $(1)$ we get the required result.
\eepf
\bl\label{re21}
Let $R=\bigoplus\limits_{\n\in\NN^s}R_{\n}$ be a standard $\NN^s$-graded Noetherian ring defined over a local ring $(R_{\bf 0},\mm)$ and $M=\bigoplus\limits_{\n\in\ZZ^s}M_{\n}$ be a finitely generated $\ZZ^s$-graded $R$-module. Let $a\in R_{\m}$ be an $M$-filter-regular where $\m\neq{\bf 0}.$ Then for all $\n\in\ZZ^s$ and $i\geq 0,$ the following sequence is exact
\beqn \displaystyle [H^{i}_{R_{++}}(M)]_{\n}\longrightarrow \left[H^{i}_{R_{++}}\lf\frac{M}{aM}\rg\right]_{\n}\longrightarrow [H^{i+1}_{R_{++}}(M)]_{\n-\m}.\eeqn 
\el
\bpf
Consider the following short exact sequence of $R$-modules $$0\longrightarrow (0:_{M}a)\longrightarrow M\longrightarrow \frac{M}{(0:_{M}a)}\longrightarrow 0.$$ Since $a$ is $M$-filter-regular, $(0:_{M}a)$ is $R_{++}$-torsion. Hence $$H^{i}_{R_{++}}(M)\simeq H^{i}_{R_{++}}\lf\frac{M}{(0:_{M}a)}\rg\mbox{ for all } i\geq 1.$$ Therefore the short exact sequence of $R$-modules $$0\longrightarrow \frac{M}{(0:_{M}a)}(-\m)\overset{a}\longrightarrow M\longrightarrow \frac{M}{aM}\longrightarrow 0$$ gives the desired exact sequence. 
\eepf
\bl\label{need}
Let $R=\bigoplus\limits_{\n\in\NN^s}R_{\n}$ be a standard $\NN^s$-graded Noetherian ring defined over a local ring $(R_{{\bf 0}},\mm)$ and $M=\bigoplus\limits_{\n\in\ZZ^s}M_{\n}$ be a finitely generated $\ZZ^s$-graded $R$-module. 
Suppose $M$ satisfies Hyry's condition $\H _R(M,{\m}).$
 Let $a_1,\dots,a_l\in R_{{\bf {e_j}}}$ be an $M$-filter-regular sequence. Then ${M}/{(a_1,\dots,a_l)M}$ satisfies Hyry's condition $\H_R({{M}/{(a_1,\dots,a_l)M},\m+l{\bf{e_j}}}).$
  
\el
\bpf
We use induction on $l.$ Let $l=1.$ By Lemma \ref{re21}, for all $i\geq 0,$ we get the exact sequence
 \beqn\label{re1}[H^{i}_{R_{++}}(M)]_{\n}\longrightarrow \left[H^{i}_{R_{++}}\lf\frac{M}{a_1M}\rg\right]_{\n}\longrightarrow [H^{i+1}_{R_{++}}(M)]_{\n-{\bf {e_j}}}.\eeqn Since for all $i\geq 0$ and $\n\geq \m+{\bf {e_j}},$ $[H^{i}_{R_{++}}(M)]_{\n}=0,$ we get $\displaystyle\left[H^{i}_{R_{++}}\lf\frac{M}{a_1M}\rg\right]_{\n}=0$ for all $i\geq 0$ and $\n\geq \m+{\bf {e_j}}.$ Hence the result is true for $l=1.$
\\ Suppose $l\geq 2$ and the result is true upto $l-1.$ Let $N=M/(a_1,\dots,a_{l-1})M.$ Then $$\displaystyle [H^{i}_{R_{++}}(N)]_{\n}=0\mbox{ for all }i\geq 0\mbox{ and }\n\geq \m+(l-1){\bf {e_j}}.$$ Since $a_l$ is $N$-filter-regular, using $l=1$ case, we get $\displaystyle \left[H^{i}_{R_{++}}\lf\frac{N}{a_lN}\rg\right]_{\n}=0$ for all $i\geq 0$ and $\n\geq \m+l{\bf {e_j}}.$   
\eepf
\bp\label{important}
Let $R=\bigoplus\limits_{\n\in\NN^s}R_{\n}$ be a standard $\NN^s$-graded Noetherian ring defined over a local ring $(R_{{\bf 0}},\mm)$ and $M=\bigoplus\limits_{\n\in\ZZ^s}M_{\n}$ be a finitely generated $\ZZ^s$-graded $R$-module.
Suppose $M$ satisfies Hyry's  condition $\H_R(M, \m).$
 Let ${\bf q}=(q_1,\dots,q_s)\in\NN^s,$ $a_{j1},\dots,a_{jq_j}\in R_{{\bf {e_j}}}$ for all $j=1,\dots,s$ such that $a_{11},\dots,a_{1q_1},\dots,a_{s1},\dots,a_{sq_s}$ be an $M$-filter-regular sequence. Then ${M}/{(a_{11},\dots,a_{1q_1},\dots,a_{s1},\dots,a_{sq_s})M}$ satisfies Hyry's condition $\H_R({{M}/{(a_{11},\dots,a_{1q_1},\dots,a_{s1},\dots,a_{sq_s})M},\m+{\bf{q}}}).$ 
\ep
\bpf
Define $N_0=M$ and $N_j={N_{j-1}}/{(a_{j1},\dots,a_{jq_j})N_{j-1}}$ for all $j=1,\dots,s.$ Since $a_{j1},\dots,a_{jq_j}\in R_{{\bf {e_j}}}$ is an $N_{j-1}$-filter-regular sequence for all $j=1,\dots,s,$ by Lemma \ref{need}, we get the required result.
\eepf
\bt\label{second}
Let $(R,\mm)$ be a Noetherian local ring of dimension $d\geq 1$  and $\idl$ be $\mm$-primary ideals in $R.$ Let $\mathcal{F}=\lbrace\fil(\n)\rbrace_{\n\in{\ZZ}^s}$ be an $\I$-admissible filtration of ideals in $R.$ 
Suppose $G (\fil)$ satisfies Hyry's condition $\H _{G (\I)}(G (\fil), \m).$
Let ${\bf q}\in\NN^s$ such that $|{\bf q}|=d$ and $\lbrace a_{ij}\in I_i:j=1,\ldots,q_i;\ii\rbrace$ be a joint reduction of $\fil$ of type ${\bf q}$ such that $a^*_{11},\dots,a^*_{1q_1},\dots,a^*_{s1},\dots,a^*_{sq_s}$ is $G(\fil)$-filter-regular sequence where $a^*_{ij}$ is the image of $a_{ij}$ in $G(\I)_{{\bf {e_i}}}$ for all $j=1,\dots,q_i$ and $\ii.$ Then $$\fil(\n)=\sum_{i=1}^{s}{\sum_{j=1}^{q_i}{a_{ij}\fil(\n-{\bf {e_i}})}}\mbox{ for all }\n\geq \m+{\bf q}.$$
\et
\bpf
By Proposition \ref{important}, we get $$\displaystyle \left[H^{i}_{G(\I)_{++}}\lf\frac{G(\fil)}{(a^*_{11},\dots,a^*_{1q_1},\dots,a^*_{s1},\dots,a^*_{sq_s})G(\fil)}\rg\right]_{\n}=0$$ for all $i\geq 0$ and $\n\geq\m+ {\bf q}.$ Since $\lbrace a_{ij}\in I_i:j=1,\ldots,q_i;\ii\rbrace$ is a joint reduction of $\fil,$ ${G(\fil)}/{(a^*_{11},\dots,a^*_{1q_1},\dots,a^*_{s1},\dots,a^*_{sq_s})G(\fil)}$ is $G(\I)_{++}$-torsion. Thus $$H^{0}_{G(\I)_{++}}\lf\frac{G(\fil)}{(a^*_{11},\dots,a^*_{1q_1},\dots,a^*_{s1},\dots,a^*_{sq_s})G(\fil)}\rg=\frac{G(\fil)}{(a^*_{11},\dots,a^*_{1q_1},\dots,a^*_{s1},\dots,a^*_{sq_s})G(\fil)}.$$ Hence for all $\n\geq\m+{\bf q},$ we have  \beqnn\label{mnb}\fil(\n)=\sum_{i=1}^{s}{\sum_{j=1}^{q_i}{a_{ij}\fil(\n-{\bf {e_i}})}}+\fil(\n+{\bf e}).\eeqnn Now for all $k\geq 0$ and  $\n\geq\m+{\bf q},$ $$\fil(\n+k{\bf e})= \sum_{i=1}^{s}{\sum_{j=1}^{q_i}{a_{ij}\fil(\n+k{\bf {e-e_i}})}}+\fil(\n+(k+1){\bf e})\subseteq\sum_{i=1}^{s}{\sum_{j=1}^{q_i}{a_{ij}\fil(\n-{\bf {e_i}})}}+\fil(\n+(k+1){\bf e}).$$
Since $\fil$ is an $\I$-admissible filtration, by \cite{rees3}, for each $\ii,$ there exists integer $r_i\in\NN$ such that for all $\n\in\NN^s,$ where $n_i\geq r_i,$ we have $\fil(\n+{\bf {e_i}})=I_i\fil(\n).$ Let $r=\max\{r_i:\ii\}.$ Therefore \beqn\fil(\n+{\bf {e}})&\subseteq&\sum_{i=1}^{s}{\sum_{j=1}^{q_i}{a_{ij}\fil(\n+{\bf e}-{\bf {e_i}})}}+\fil(\n+2{\bf {e}})\subseteq\sum_{i=1}^{s}{\sum_{j=1}^{q_i}{a_{ij}\fil(\n-{\bf {e_i}})}}+\fil(\n+2{\bf {e}})\\&\subseteq\vdots&\\&\subseteq &\sum_{i=1}^{s}{\sum_{j=1}^{q_i}{a_{ij}\fil(\n-{\bf {e_i}})}}+\fil(\n+(r+2){\bf {e}})\subseteq \sum_{i=1}^{s}{\sum_{j=1}^{q_i}{a_{ij}\fil(\n-{\bf {e_i}})}}+\id\fil(\n).\eeqn Hence by using Nakayama's Lemma, from the equality (\ref{mnb}) we get the required result.
\eepf
\bl\label{vanishing}
Let $(R,\mm)$ be a Noetherian local ring of dimension $d\geq 1$  and $\idl$ be $\mm$-primary ideals in $R.$ Let $\mathcal{F}=\lbrace\fil(\n)\rbrace_{\n\in{\ZZ}^s}$ be an $\I$-admissible filtration of ideals in $R.$ Suppose $\R (\fil)$ satisfies Hyry's condition $\H_{\R{(\I)}}(\R(\fil),\m)$ where $\m\in\NN^s.$ Then $G (\fil)$ satisfies Hyry's condition $\H _{G (\I)}(G (\fil), \m).$
\el
\bpf
Denote $\mathcal{R}'(\fil)/\mathcal{R}'(\fil)({\bf e})$ by $G'(\fil).$ Consider the short exact sequence of $\R(\I)$-modules 
 \begin{equation}\label{equation}0\longrightarrow \mathcal{R}'(\fil)({\bf e})\longrightarrow \mathcal{R}'(\fil)\longrightarrow G'(\fil)\longrightarrow 0.
 \end{equation} This induces the long exact sequence of $R$-modules
$$\cdots\longrightarrow [H_{{\R}_{++}}^i(\mathcal{R}'(\fil))]_{\n+{\bf e}} \longrightarrow [H_{{\R}_{++}}^i(\mathcal{R}'(\fil))]_{\n} 
\longrightarrow [H_{{\R}_{++}}^i(G'(\fil))]_{\n}
 \longrightarrow [H_{{\R}_{++}}^{i+1}(\mathcal{R}'(\fil))]_{\n+{\bf e}} \longrightarrow \cdots.$$ Hence by \cite[Proposition 4.2]{msv} and change of ring principle, we get the required result.  
\eepf
\bt\label{final}
Let $(R,\mm)$ be a Noetherian local ring of dimension $d\geq 1$ and $\idl$ be $\mm$-primary ideals in $R.$ Let $\mathcal{F}=\lbrace\fil(\n)\rbrace_{\n\in{\ZZ}^s}$ be an $\I$-admissible filtration of ideals in $R$ and $\R(\fil)$ satisfy Hyry's condition  $\H_{\R{(\I)}}(\R(\fil), {\bf 0}).$ Let ${\bf q}\in\NN^s$ such that $|{\bf q}|=d.$ Then there exists a joint reduction $\lbrace a_{ij}\in I_i:j=1,\ldots,q_i;\ii\rbrace$ of $\fil$ of type ${\bf q}$ such that $$\fil(\n)=\sum_{i=1}^{s}{\sum_{j=1}^{q_i}{a_{ij}\fil(\n-{\bf {e_i}})}}\mbox{ for all }\n\geq {\bf q}\mbox{ ~~~~and }$$ $$\jr_{{\bf q}}(\fil)\leq\max\{q_i\mid i\in A\}-1\mbox{ ~where~ } A=\{i\mid q_i\geq 1\}.$$
\et
\bpf
Since $\R(\fil)$ satisfies Hyry's condition  $\H_{\R{(\I)}}(\R(\fil),{\bf 0}),$ by Lemma \ref{vanishing}, $G(\fil)$ satisfies Hyry's condition  $\H_{G(\I)}(G(\fil), {\bf 0}).$ Therefore by Theorem \ref{joint}, there exists a joint reduction $\lbrace a_{ij}\in I_i:j=1,\ldots,q_i;\ii\rbrace$ of $\fil$ of type ${\bf q}$ such that such that 
$a^*_{11},\dots,a^*_{1q_1},\dots,a^*_{s1},\dots,a^*_{sq_s}$ is a $G(\fil)$-filter-regular sequence where $a^*_{ij}$ is the image of $a_{ij}$ in $G(\I)_{\bf {e_i}}$ for all $j=1,\ldots,q_i,$ $\ii.$ Hence by Theorem \ref{second},  for all $\n\geq {\bf q},$  $$\fil(\n)=\sum\limits_{i=1}^{s}{\sum\limits_{j=1}^{q_i}{a_{ij}\fil(\n-{\bf {e_i}})}}.$$ 
\eepf
\bex
\rm{Let $R=k[|X,Y|].$ Then $R$ is a regular local ring of dimension two. Let $I=(X,Y^2)$ and $J=(X^2,Y).$ Then $I,J$ are complete parameter ideals in $R.$ Consider the filtration $\fil=\{I^rJ^s\}_{r,s\in\ZZ}.$ Since $I,J$ are complete ideals, by \cite[Theorem 2$'$, Appendix 5]{zariski}, $I^r,$ $J^s$ and $I^rJ^s$ are complete ideals for all $r,s\geq 1.$ By \cite[Proposition 3.2]{msv} and \cite[Proposition 3.5]{msv}, for all $r,s\in\NN,$ $H_{\R_{++}}^1({\R}(I,J))_{(r,s)}=0.$ Note that $(X^3+Y^3,XY)$ is a minimal reduction of $IJ$ and $$(X^3+Y^3,XY)IJ=I^2J^2.$$ Thus $r(IJ)\leq 1$ and hence $e_2(IJ)=0.$ Therefore using \cite[Theorem 4.3]{msv} and \cite[Lemma 2.11]{msv}, we get $[H_{\R_{++}}^2({\R}(I,J))]_{(r,s)}=0$ for all $r,s\in\NN.$ Hence  $\R(\fil)$ satisfies Hyry's condition $\H _{\R{(\I)}}(\R(\fil), {\bf 0}).$
\\ Note that $\mathcal A=\{X,Y\}$ is a joint reduction of $(I,J)$ of type ${\bf e}$ and 
$$XI^rJ^{s+1}+YI^{r+1}J^s=I^{r+1}J^{s+1}\mbox{ for all } r,s\in\NN.$$ 
Thus $\jr_{\bf e}(\fil)=\jr_{\mathcal A_{\bf e}}(\fil)=0.$
}
\eex
\bt\label{complete}
Let $(R,\mm)$ be an analytically unramified Noetherian local ring of dimension $d\geq 2$  and let $\idl$ be $\mm$-primary ideals in $R.$ Let $\overline{\R}(\I)=\bigoplus\limits_{\n\in\NN^s}\overline{\I^{\n}}$ satisfy Hyry's condition $\H_{\R{(\I)}}(\overline{\R}(\I),{\bf 0}).$ Suppose $\I^{\n}$ is complete for all $\n\in\NN^s$ such that $1\leq|\n|\leq d-1.$ Then $\I^{\n}$ is complete for all $\n\in\NN^s$ with $|\n | \geq 1.$ 
\et
\bpf We use induction on $|\n|.$ By given hypothesis the result is true upto $1\leq|\n|\leq d-1.$ Suppose $\n\in\NN^s$ with $|\n|\geq d$ and the result is true for all ${\bf k}\in\NN^s$ such that $1\leq|{\bf k}|<|\n|.$ Let $\m=(m_1,\dots,m_s)\in\NN^s$ such that $\m\leq\n$ and $|\m|=d.$ Consider the filtration $\fil=\lbrace \overline{\I^{\n}} \rbrace_
{\n\in\ZZ}.$ By \cite{rees3}, $\fil$ is an $\I$-admissible filtration. Then by Theorems \ref{joint} and \ref{second}, there exists a joint reduction $\lbrace a_{ij}\in I_i:j=1,\ldots,m_i;\ii\rbrace$ of $\fil$ of type $\m$
such that $$\overline{\I^{\rrr}}=\sum_{i=1}^{s}{\sum_{j=1}^{m_i}{a_{ij}\overline{\I^{\rrr-{\bf{e_i}}}}}}\mbox{ for all } \rrr\geq\m.$$ Thus $\overline{\I^{\n}} =\sum\limits_{i=1}^{s}{\sum\limits_{j=1}^{m_i}{a_{ij}\overline{\I^{\n-{\bf{e_i}}}}}}.$  By induction hypothesis, $\I^{\n-{\bf{e_i}}}$ is complete for all $i\in A: =\{i | n_i \geq 1\}.$ Hence $$\overline{\I^{\n}} =\sum_{i=1}^{s}{\sum_{j=1}^{m_i}{a_{ij}\overline{\I^{\n-{\bf {e_i}}}}}}=\sum_{i=1}^{s}{\sum_{j=1}^{m_i}{a_{ij}{\I^{\n-{\bf {e_i}}}}}}\subseteq \I^{\n}.$$
\eepf
As a consequence of the above theorem we obtain a generalisation of a theorem of Reid, Roberts and Vitulli \cite[Proposition 3.1]{V} about complete monomial ideals.
\bt
Let $R=k[X_1,\dots,X_d]$ where $d\geq 1$ and $\mm=(X_1,\dots,X_d)$ be the maximal homogeneous ideal of $R.$ Let $\idl$ be $\mm$-primary monomial ideals of $R.$ Suppose $\I^{\n}$ is complete for all $\n\in\NN^s$ such that $1\leq|\n|\leq d-1.$ Then $\I^{\n}$ is complete for all $\n\in\NN^s$ with $|\n | \geq 1.$
\et
\bpf
If $d=1$ then $R$ is a PID and hence normal. Therefore every ideal is complete since principal
ideals in normal domains are complete. Let $d \geq 2.$ 
Since $\idl$ are monomial ideals, $\overline{\R}(\I)$ is Cohen-Macaulay by \cite[Theorem 6.3.5]{BH}. Let $W=R\setminus\mm.$ Then $S=W^{-1}\overline{\R}(\I)$ is Cohen-Macaulay. Since for any ideal $I,$ $W^{-1}\overline I=\overline{W^{-1}I},$ we have $$W^{-1}\overline{\R}(\I)=\bigoplus\limits_{\n\in\NN^s}W^{-1}\overline{\I^{\n}}=\bigoplus\limits_{\n\in\NN^s}(\overline{W^{-1}(\I^{\n})})=\overline{\R}(W^{-1}I_1,\dots,W^{-1}I_s).$$ Therefore $S$ satisfies Hyry's condition $\H_Q(S,{\bf 0})$ where $Q=W^{-1}{\R}(\I).$  Replace $R$ by $W^{-1}R.$ Therefore by Theorem \ref{complete}, $W^{-1}(\I^{\n})$ is complete for all $\n\in\NN^s$ such that $|\n|\geq 1.$ Since $\mm$ is the maximal homogeneous ideal of $R$ and $\displaystyle W^{-1}\lf{\overline{\I^{\n}}}/{{\I^{\n}}}\rg=0,$ we get the required result.
\eepf

We end the paper with three examples illustrating some of the results proved above.
\bex
{\rm Let $S=\QQ[[X,Y,Z]],$ $f=X^2+Y^2+Z^2.$ Then $R=S/(f)$ is analytically unramified Cohen-Macaulay   reduced local ring of dimension $2.$ Set $\mm=(X,Y,Z)/(f).$ Since $G_{\mm}(R)=\bigoplus\limits_{n\geq 0}\mm^n/\mm^{n+1}\simeq 
\QQ[X,Y,Z]/(f)$ is reduced,  $\mm^n$ is complete for all $n\geq 1.$
\\ Consider the  $\mm$-admissible  filtration $\fil=\{\ov{\mm^n}\}_{n\in\ZZ}.$  The Hilbert polynomial of the filtration $\fil$ is $\po(n)=2\binom{n+1}{2}-n.$ 
Set $\mathcal R=\mathcal R(\mm).$ 
Since $R$ is Cohen-Macaulay, $H^0_{\mathcal R_{++}}({\R}(\fil))=0.$ By \cite[Theorem 3.5]{blancafort}, $[H^1_{\mathcal R_{++}}({\R}(\fil))]_n=\widetilde{\ov{\mm^n}}/{\ov{\mm^n}}$ for all $n\geq 0$ where $\{\widetilde{\ov{\mm^n}}\}_{n\in\ZZ}$ is the Ratliff-Rush closure filtration of $\fil.$ Therefore by \cite[Proposition 3.2]{msv}, 
$[H^1_{\mathcal R_{++}}({\R}(\fil))]_n=0$ for all $n\geq 0.$ By \cite[Theorem 4.1]{blancafort}, we get $[H^2_{\mathcal R_{++}}({\R}(\fil))]_0=0.$ Hence by \cite[Lemma 4.7]{blancafort}, $[H^2_{\mathcal R_{++}}({\R}(\fil))]_n=0$ for all $n\geq 0.$ Hence ${\R}(\fil)$ satisfies the condition $\H_{\R(\mm)}({\R}(\fil),0).$
}\eex
The following examples show that Hyry's  condition $\H_{\R{(\I)}}(\overline{\R}(\I),{\bf 0})$ is sufficient but  not necessary in Theorem \ref{complete}. 
\bex{\rm
Let $S=\QQ[[X,Y,Z]],$ $g=X^3+Y^3+Z^3.$ Then $R=S/(g)$ is analytically unramified Cohen-Macaulay reduced local ring of dimension $2.$ Set  $\mm=(X,Y,Z)/(g).$  Since  $G_{\mm}(R)=\bigoplus\limits_{n\geq 0}\mm^n/\mm^{n+1}\simeq \QQ[X,Y,Z]/(g)$ is reduced,  $\mm^n$ is complete for all $n\geq 1.$
\\ Consider the $\mm$-admissible filtration $\fil=\{\ov{\mm^n}\}_{n\in\ZZ}.$ The Hilbert polynomial of the filtration $\fil$ is $\po(n)=3\binom{n+1}{2}-3n+1.$ 
Set $\mathcal R=\mathcal R(\mm).$ 
Since $R$ is Cohen-Macaulay, $H^0_{\mathcal R_{++}}({\R}(\fil))=0.$ By \cite[Theorem 3.5]{blancafort}, $[H^1_{\mathcal R_{++}}({\R}(\fil))]_n=\widetilde{\ov{\mm^n}}/{\ov{\mm^n}}$ for all $n\geq 0$ (here $\{\widetilde{\ov{\mm^n}}\}_{n\in\ZZ}$ is the Ratliff-Rush closure filtration of $\fil$). Therefore $[H^1_{\mathcal R_{++}}({\R}(\fil))]_0=0$ for all $n\geq 0.$ By \cite[Theorem 4.1]{blancafort}, we get $\lm\lf H^2_{\mathcal R_{++}}({\R}(\fil))\rg_0=1.$ Hence ${\R}(\fil)$ does not satisfy Hyry's condition $\H_{\R(\mm)}({\R}(\fil),0).$
}\eex
\bex{\rm
Let $R=k[x,y,z]$ where $k$ is a field of characteristic not equal to $3$ and $I=(x^4,x(y^3+z^3),y(y^3+z^3),z(y^3+z^3))+(x,y,z)^5.$ In \cite[Theorem 3.12]{HH}, S. Huckaba and Huneke showed that $\height(I)=3,$ $I$ is normal ideal, i.e. $\ov{I^n}=I^n$ for all $n\geq 1$ and $H^2(X,\mathscr{O}_X)\neq 0$ where $X=\Proj \R(I).$ Hence $H^3_{\mathcal R_{++}}({\R}(I))_0\neq 0.$ Thus ${\R}(I)$ does not satisfy Hyry's condition $\H_{\R(I)}({\R}(I),0).$
}\eex

\end{document}